\documentclass{amsart}

\usepackage{amsmath,
             amssymb, amsfonts}

\usepackage{graphicx}

\newtheorem{theorem}{Theorem}[section]
\newtheorem{corollary}[theorem]{Corollary}
\newtheorem{lemma}[theorem]{Lemma}
\newtheorem{proposition}[theorem]{Proposition}

\theoremstyle{definition}

\theoremstyle{remark}
\newtheorem{remark}[theorem]{Remark}

\numberwithin{equation}{section}

\newcommand{\norm}[1]{\left\lVert#1\right\rVert}%

\newcommand{\rom}[1]{%
  \textup{\uppercase\expandafter{\romannumeral#1}}%
}

\begin{document}

\title[Euler equations
 in the end-point critical Triebel-Lizorkin space]
 {Persistence of the solution to the Euler equations
 in the end-point critical Triebel-Lizorkin space
$F^{d+1}_{1, \infty}(\mathbb{R}^d)$ }


\author{Hee Chul Pak}
\address{Department of Mathematics,
         Dankook University, 119 Dandae-ro, Dongnam-gu, Cheonan-si, Chungnam, 31116,  Republic of Korea}
\email{hpak@dankook.ac.kr}
\thanks{Correspondence: Hee Chul Pak, hpak@dankook.ac.kr}


\author{Jun Seok Hwang}
\address{Department of Mathematics,
         Dankook University, 119 Dandae-ro, Dongnam-gu, Cheonan-si, Chungnam, 31116,  Republic of Korea}




\begin{abstract}
Local stay of
 the solutions to the  Euler equations for an ideal incompressible
  fluid
in the end-point  Triebel-Lizorkin spaces $F^s_{1, \infty}(\mathbb{R}^d)$
with  $s \geq d + 1$ is clarified.
\end{abstract}

\maketitle

\section{Introduction}

The perfect  incompressible inviscid fluid is governed by the Euler
equations:
\begin{align}
 \frac{\partial}{\partial t} \; u
  + (u, \nabla) u   &= - \nabla p
                                      \label{Euler}\\
\mbox{div }u
     & = 0. \label{Euler-2}
\end{align}
Here $u(x,t) = (u_1, u_2, \cdots , u_d)$ is the 
{\it velocity} of a fluid flow and $p(x, t)$ is the scalar {\it
pressure}.

Existence and uniqueness theories  of
solutions of  the 2 or 3 dimensional Euler equations
have been worked on by many mathematicians and  physicists.
L. Lichtenstein, N. Gunther and
Wolibner started the subject on H\"{o}lder classes.
 D. Ebin and J. Marsden, J. Bourguignon and
H. Brezis, R. Temam, T. Kato and G. Ponce studied this subject on
Sobolev spaces.
Some researches on the Euler equations
in  Besov spaces $B_{p, q}^{s}(\mathbb{R}^d)$  have been done by M.
Vishik \cite{V1}, D. Chae \cite{Chae2}  and Pak-Park \cite{P-P1}.
For a detailed survey of this
issue, we refer \cite{B-C-D, Majda-Bertozzi, C, Constantin, Bourgain-Li-1, Bourgain-Li-2}
and references therein.
Bourgain and Li proved strong ill-posedness results for
the Euler equations associated with initial data in (borderline) Besov spaces,
Sobolev spaces or the space $C^m$.
For the survey of the ill-posedness issue we refer
\cite{Bourgain-Li-1, Bourgain-Li-2}.

The existence and uniqueness of the solutions of the Euler equations
in general Triebel-Lizorkin spaces $F_{p, q}^{s}(\mathbb{R}^d)$ were
first investigated by D. Chae in \cite{Chae1}, and his research
covers the cases of
\[
s > \frac{d}{p} + 1 , \; 1 < p, q < \infty \quad \mbox{ or }  \quad
 s \geq d + 1, \; p= 1,  1 < q < \infty.
\]
The unique existence of solution in the critical space
$F_{1, 1}^{s}(\mathbb{R}^d) = B_{1, 1}^{s}(\mathbb{R}^d)$
($s \geq d + 1$)
is proved by Pak-Park in \cite{P-P3}, and very recently,
the continuity of the solution map
in this space is reported in \cite{Guo-Li}. The
remaining critical and sub-critical cases
\[
s = \frac{d}{p} + 1 , \; 1 < p \leq \infty, q=1 \quad \mbox{ and }
\quad
 s \geq d + 1, \; p= 1,  q = \infty
\]
of well- or ill-posedness issue are not clearly reported, yet.
In this paper, we clarify the persistence of
the solution
for the critical end-point case
 $p = 1$ and $q = \infty$ with $s \geq d + 1$.
We now state our results as follows.
\begin{theorem} \label{Main Theorem}
Let $s$ be a real number greater than or equal to $d +1$.   \\
1. (Persistence and uniqueness)
For any divergence free vector field
$u_0 \in  F_{1,\infty}^{s}(\mathbb{R}^d)$, there exists a positive time
$T > 0$ such that
 the initial value problem for the
Euler equations (\ref{Euler}) with initial velocity $u_0$ has  a
unique solution $u$
with
$u(t) \in F_{1,\infty}^{s}(\mathbb{R}^d)$ for all $t \in [0, T]$.
\\
Moreover, this solution $u$ also satisfies    \\
2. (2-D global existence of the solution) When $d =2$, the solution
$u$ does not blow-up within finite time.
That is,
$u(t) \in
F_{1,\infty}^{s}(\mathbb{R}^2) $
for all $t \in [0, \infty)$.
\end{theorem}

One of the main analytic  points to deal with the end-point
critical spaces is
that the usual techniques may not be utilized in those critical spaces.
For example,
the (vector-valued) Hardy-Littlewood maximal operator
is not $L^1$-bounded and even  the Calderon-Zygmund theory does not work.
Furthermore,
 the smooth approximation
via mollifier  can not be applicable
due to the fact that  the Schwartz class is not dense
in the $\ell^{\infty}$-hierarchy spaces, which  causes
uncertainty of
 the existence result in the space $ F_{1,\infty}^{s}(\mathbb{R}^d)$.
However,
we  could detour and organize some known estimates and  techniques
to derive our conclusion without too many analytic difficulties.

We pursue  self-contained
and detailed descriptions of the arguments to avoid possible gaps
that sometimes happen on this subject\cite{Ch, Guo-Li}.

\bigskip

\noindent{\bf Notations:}  \label{notation} Throughout this paper,
\begin{enumerate}
\item
for $x \in  \mathbb{R}^d$, $x_i$ is the $i$-th component of $x$
\item
$B(x, r) := \{ y \in \mathbb{R}^d  :   |x- y| < r  \}$
\item
$\displaystyle  \frac{\partial f}{\partial x_{k}} =
\partial_{x_k} f$ or simply $\partial_{k} f$
\item
 for $k > 0$ and a function $\phi$ on $\mathbb{R}^d$,
$[\phi]_{k} (x) := k^{d} \phi(k x)$ for $x \in \mathbb{R}^d$
\item
for $f \in \mathcal{S}({\mathbb R}^{d})$, the Fourier transform
$\hat{f}= {\mathcal F}(f)$ of $f$ on ${\mathbb R}^{d}$ is defined by
\[
\hat{f}(\xi) = {\mathcal F}(f)(\xi) =  \int_{{\mathbb R}^{n}} f(x)
e^{- i x \cdot \xi}  \, dx
\]
\item
$D^s u = \sqrt{ -\Delta}^{\;s} u :=  {\mathcal F}^{-1} (| \cdot |^s
{\mathcal F}(u)) $
\item
for $1 \leq p < \infty$
and a sequence of functions
$
\left\{
              f_j
\right\}\!\!_{j \in \mathbb{Z}}$,
\[
        \left\|
\left\{
              f_j
\right\}_{j \in \mathbb{Z}}
        \right\|_{l^p}\!\!(x)
        :=
        \left(
        \sum_{j \in \mathbb{Z}}  \left|f_j(x)  \right|^p
        \right)^{1/p}
\]
\item
the notation $X \lesssim Y$ means that $X \leq C  Y$, where $C$ is a
fixed but unspecified constant. Unless explicitly stated otherwise,
$C$ may depend on the dimension $d$ and various other parameters
(such as exponents), but not on the functions or variables $(u, v,
f, g, x_i, \cdots)$ involved.

\end{enumerate}

\section{A-priori estimate}
\label{prelim}

Let $\mathcal{S}(\mathbb{R}^d)$ denote the Schwartz class. We consider a
nonnegative radial function $\chi \in \mathcal{S}(\mathbb{R}^d)$
satisfying
$ \mbox{supp} \:\chi \subset \{ \xi \in {\mathbb{R}}^d :
|\xi| \leq 1 \}$, and $\chi =1$
for $|\xi| \leq \frac{3}{4}$. \label{chai-1}
Set $h_j(\xi) := \chi(2^{-j-1}\xi) - \chi( 2^{-j}
\xi)$
and let $\varphi_j$ and $\Phi$ be defined by $\varphi_j:=
\mathcal{F}^{-1}(h_j)$, $j \geq 0$ and
$\Phi :=
\mathcal{F}^{-1}(\chi)$.
 For any $f\in\mathcal{S}'(\mathbb{R}^{d})$, we define the
operator $\Delta_{j}$ by
\begin{align*}
\Delta_{j}f=
\begin{cases}
\hat{\varphi_{j}}(D)f=\varphi_{j}\ast f    &\text{for } j \geq 0 \\
\hat{\Phi}(D)f=\Phi\ast f                 &\text{for } j = -1
\\
0                                          &\text{for } j \leq -2
\end{cases}
\end{align*}
and the operator $\dot{\Delta}_{j}$ by
\begin{align*}
\dot{\Delta}_{j} f =
\begin{cases}
\Delta_{j}f                               &\text{for } j \geq 0 \\
\hat{\varphi_{j}}(D)f=\varphi_{j} \ast f                  &\text{for
} j \leq -1.
\end{cases}
\end{align*}
Then we have an analog of a partition of unity:
\begin{align*}
\Delta_{-1} + \sum_{j = 0}^{\infty} \Delta_j = I \quad \mbox{ and  }
\quad \sum_{j = - \infty}^{\infty} \dot{\Delta}_j = I,
\label{part-1-4}
\end{align*}
where $I$ represents the identity operator.
The partial sum operators $S_k f$ and $\dot{S}_k f$ are
$\displaystyle{S_{k}f=\sum_{j=-\infty}^{k}\Delta_{j}f}$
and $\displaystyle{\dot{S}_{k}f=\sum_{j=-\infty}^{k} \dot{\Delta}_{j}f}$,
respectively.

For $s\in\mathbb{R}$,
the homogenous Triebel-Lizorkin space
$\dot{F}_{1,\infty}^{s}(\mathbb{R}^d)$ is the collection of all $f
\in \mathcal{S}'(\mathbb{R}^d)$ modulo polynomials such that
\[
\norm{f}
  _{\dot{F}_{1,\infty}^{s}}
:= \int
  _{\mathbb{R}^{d}}
\underset{j\in\mathbb{Z}}
  {\sup}
     \left|
          2^{js} \dot{\Delta}_{j}f
     \right|
(x)dx<\infty,
\]
and the nonhomogeneous Triebel-Lizorkin space
$F_{1,\infty}^{s}(\mathbb{R}^d)$ is the space of all tempered
distributions $f \in \mathcal{S}'(\mathbb{R}^d)$ obeying
\begin{align}
\|f \|
  _{F_{1,\infty}^{s}}
:=
\int
  _{\mathbb{R}^{d}}
\underset{j\in\mathbb{Z}}
  {\sup}
     \left|
          2^{js} {\Delta}_{j}f
     \right|
(x) dx < \infty.   \label{t-norm-1}
\end{align}
We observe that for $s > 0$, the Triebel-Lizorkin norm
$\norm{f}_{F_{1,\infty}^{s}}$ is equivalent to the nonhomogeneous
norms
\begin{equation}
  \norm{f}_{L^1}
  +\norm{f}_{\dot{F}_{1,\infty}^{s}}.
\end{equation}

\begin{remark} \label{rmk-2}
Let $1 \leq p < \infty$, $1 \leq q \leq \infty$ and $s \in
\mathbb{R}$. For any $k \in \mathbb{N}$, we have Bernstein type
inequality:
\[
\norm{ D^k  f }_{\dot{F}^{s}_{p, q}} \lesssim \norm{  f
}_{\dot{F}^{s+k}_{p, q}} \lesssim \norm{ D^k  f }_{\dot{F}^{s}_{p,
q}}.
\]
\end{remark}


\subsection{Basic estimates}
\label{Basic estimates}

We begin with Lemma \ref{Peetre} which is introduced by Peetre.
Lemma \ref{Peetre} has been developed for the
proof of the independent choice of the mother (bump) function $\chi$
employed at page \pageref{prelim}. The version given in this paper
is a refinement introduced by Triebel \cite{T}. For the proof, we
refer the page 71 in \cite{Grafakos-modern}.

\begin{lemma}  \label{Peetre}
Let $0<  r <\infty$. Then
for all $t>0$ and for any $C^{1}$-function $u$ on $\mathbb{R}^{d}$
whose Fourier transform is supported in the ball $|\xi| \lesssim t$
and that satisfies
\begin{align}
|u(z)| \lesssim (1+|z|)^{\frac{d}{r}}, \label{struc-cond}
\end{align}
 we have
\begin{equation}\label{eq6.5.6}
\underset{z\in\mathbb{R}^{d}}{\sup}
    \frac{1}
         {t}
    \frac{|
         \nabla u(x-z)|}
         { (1 + t|z|)^\frac{d}{r}}
\lesssim
    \underset{z\in\mathbb{R}^{d}}{\sup}
    \frac{|u(x-z)|}
         { (1+ t|z|)^\frac{d}{r}}
\lesssim
    M(|u|^{r})(x)^{\frac{1}{r}},
\end{equation}
where $M$ is the Hardy-Littlewood maximal operator.
\end{lemma}

Paley-Wiener-Schwartz theorem says that every
distribution $u$ on $\mathbb{R}^{d}$ whose Fourier transform is
supported in the ball $|\xi|\leq c_{0}t$ is an entire function of
$d$-complex variables.
The following proposition is a
generalizations of one presented by Guo and Li \cite{Guo-Li}.

\begin{proposition} \label{prop:3}
Let $r, r_1, r_2  > 0$, $s, t, \ell > 0$ with $\frac{s}{t}< \ell$
and $\gamma, \delta > 0$ with $\gamma+\delta=1$ and
 $\frac{\gamma}{r_1}  + \frac{\delta}{r_2} = \frac{1}{r}$.
Let $\psi$ be a measurable function such that
\[
|\psi(z)|
    \left(
        1 + |z|
    \right)^{d/r}
\leq U(z), \qquad z\in\mathbb{R}^d
\]
for some nonnegative radial decreasing integrable function $U$. Then
we have that for $x\in \mathbb{R}^d$,
\begin{equation}
\left|
      [\psi]_s \ast
      \left(
        gf
      \right)
\right|(x) \lesssim \left(
    \frac{t}{s}
\right)^{\frac{d}{r}} M(g)(x) \left[
    M(|f|
        ^{r_1   })
    (x)
\right]^{ \frac{\gamma}{r_1} } \left[
    M(|f|
        ^{ r_2    })
    (x)
\right]^{ \frac{\delta}{r_2} }
\end{equation}
for all $g \in L^{1}_{loc}(\mathbb{R}^d)$ and any tempered
distribution $f\in\mathcal{S}'(\mathbb{R}^d)$ whose frequency support
$\mbox{supp} \: \widehat{f}$ is contained in the ball $|\xi| \lesssim t$ together
with the condition $ |f(x)|
    \lesssim
    (1+|x|)^{d \theta}
$ ($\theta = \min\{\frac{\gamma}{r_1}, \frac{\delta}{r_2} \}$) for
all $x \in \mathbb{R}^d$.
\end{proposition}

\noindent {\bf Proof.} We take $g \in L^{1}_{loc}(\mathbb{R}^d)$ and
a distribution $f$ satisfying $\mbox{supp} \: \widehat{f} \,
  \subseteq B(0, c_0 t)$
for some $c_0>0$ and $ |f(x)|
    \lesssim
    (1+|x|)^{\frac{d}{r}}
$ for all $x \in \mathbb{R}^d$. Then we have that for $x \in
\mathbb{R}^d$,
\begin{align*}
\left|
      [\psi]_s \ast
      \left(
      gf
      \right)
\right|(x) &\leq \int_{\mathbb{R}^d}
   \left|
      \psi(z)
   \right|
   \left|
     \left(
     gf
      \right)
        \left(x - s^{-1} z\right)
   \right|   dz
\\&
\leq
  \left[
    \int_{\mathbb{R}^d}
      \left|
        \psi(z)
      \right|
      \left|
      g (x - s^{-1} z)
      \right|
      \left(
      1+\frac{t}{s}|z|
      \right)^{d/r}
       dz
  \right]
\\&\qquad
\times \underset{z\in\mathbb{R}^d}{\sup} \frac{\left|
          f \left(x-  s^{-1}z\right)
      \right|}
      {\left( 1 +  t/s \left|z\right|\right)^{d/r}}.
\end{align*}
 Lemma \ref{Peetre} yields
\begin{align*}
\underset{z\in\mathbb{R}^d}{\sup} \frac{\left|
          f \left(x-  s^{-1}z\right)
      \right|}
      {\left( 1 +  t/s \left|z\right|\right)^{d/r}}
&\leq    
\left( \underset{z\in\mathbb{R}^d}{\sup} \frac{\left|
          f \left(x- s^{-1}z\right)
      \right|}
      {\left( 1 +  \frac{t}{s} \left|z\right|\right)^{\frac{d}{r_1}}}
\right)^{\gamma} \left( \underset{z\in\mathbb{R}^d}{\sup}
\frac{\left|
          f \left(x- s^{-1}z\right)
      \right|}
      {\left( 1 + \frac{t}{s} \left|z\right|\right)^{\frac{d}{r_2}}}
\right)^{\delta}
\\&
\lesssim        
\left[
    M(|f|
        ^{r_1   })
    (x)
\right]^{ \frac{\gamma}{r_1} } \left[
    M(|f|
        ^{ r_2    })
    (x)
\right]^{ \frac{\delta}{r_2} }.
\end{align*}
By the assumption that $\frac{s}{t}< \ell$,  we also get
\begin{align*}
\int_{\mathbb{R}^d}
    \left|
        \psi(z)
    \right|
&\left|
        g (x - s^{-1} z)
    \right|
    \left(
      1+\frac{t}{s}|z|
    \right)^{d/r}
 dz
\\&
   \leq
    \max \{ \ell^{- \frac{d}{r}}, 1\}
   \int_{\mathbb{R}^d}
      \left|
         \psi(z)
      \right|
      \left(
        \frac{t}{s}
        +
           \frac{t}{s}
          \left|
            z
          \right|
      \right)^{d/r}
      \times
      | g (x - s^{-1} z) |
       dz
\\&
\lesssim
    \left( \frac{t}{s}  \right)^{\frac{d}{r}}
        \left(
           [\phi]_s \ast |g|
        \right)
       (x)
   \lesssim
    \left( \frac{t}{s}  \right)^{\frac{d}{r}}
    M(g)(x),
\end{align*}
where
$ \phi(z) := |\psi(z)|
      \left(
        \frac{t}{s}
        +
           \frac{t}{s}
          \left|
            z
          \right|
      \right)^{d/r}
$
and $[\phi]_{s}(s) := s^{d} \phi  \left(s x\right)$.
\hfill$\Box$\par

\bigskip

As an application, we have the following statement.

\begin{corollary}
\label{coro:1}
Let $\ell_0, r, r_1, r_2 >0$ and  $\gamma, \delta >0$
with $\gamma+\delta=1$ and
 $\frac{\gamma}{r_1}  + \frac{\delta}{r_2} = \frac{1}{r}$.
Then for $j,\,k\in \mathbb{Z}$ with $j>k- \ell_0$, we have that
\begin{equation}
    \left|
      \Delta_k
      \left(
        gf
      \right)
\right|(x) \lesssim 2^{(j-k)\frac{d}{r}} M(g)(x) \left[
    M(|f|
        ^{r_1   })
    (x)
\right]^{ \frac{\gamma}{r_1} } \left[
    M(|f|
        ^{ r_2    })
    (x)
\right]^{ \frac{\delta}{r_2} }
\end{equation}
for any $g \in L^{1}_{loc}(\mathbb{R}^d)$ and any measurable
function $f$ whose frequency support
$ \mbox{supp} \: \widehat{f}$ is contained in the ball
$|\xi| \lesssim 2^j$ together with the condition $ |f(x)|
    \lesssim
    (1+|x|)^{d \theta}
$ ($\theta = \min\{\frac{\gamma}{r_1}, \frac{\delta}{r_2} \}$) for
all $x \in \mathbb{R}^d$.
\end{corollary}

One of main tools for the estimates of Triebel-Lizorkin spaces is
the continuity of the Hardy-Littlewood maximal operator $M$ for the
vector valued functions. In fact,
  for any sequence
$\left\{
   g_k
 \right\}_{k \in \mathbb{Z}}$ in $L^p(\mathbb{R}^d)$,
 one has
\begin{align}
\left\| \left\|
 \left\{
   M (g_k)
 \right\}_{k \in \mathbb{Z}}
\right\|_{\ell^q} \right\|_{L^p} \lesssim \left\| \left\|
 \left\{
   g_k
 \right\}_{k \in \mathbb{Z}}
\right\|_{\ell^q} \right\|_{L^p}. \label{cont_max-fn}
\end{align}
This holds for the case of $1 < p < \infty$ and $1 < q \leq
\infty$, and unfortunately not for $p=1$.
We present a useful alternative of the estimate
(\ref{cont_max-fn}) for $L^1$-hierarchy spaces:
\begin{corollary}\label{coro-2}
Let $1 < q \leq \infty$ and let $\psi$  be a measurable  function
such that $\hat{\psi}$ has a compact support and
$$
\left| \psi(z) \right|
      \left(
        1 +  |z|
      \right)^{d s  }
\leq U(z), \qquad z \in \mathbb{R}^d
$$
for some $s > 0$ and some nonnegative radial decreasing integrable
function $U$.
Then for a fixed integer $\ell$,
\[
\left\| \left\|
 \left\{
    [\psi]_{2^{k + \ell}} * \dot{\Delta}_k f_k)
 \right\}_{k \in \mathbb{Z}}
\right\|_{\ell^q} \right\|_{L^1} \lesssim \left\| \left\|
 \left\{
   \dot{\Delta}_k f_k
 \right\}_{k \in \mathbb{Z}}
\right\|_{\ell^q} \right\|_{L^1}.
\]
\end{corollary}
\noindent {\bf Proof.}
We apply Proposition \ref{prop:3} to have
\begin{align*}
&\left\| \left\|
 \left\{
    [\psi]_{2^k} * \dot{\Delta}_k f_k)
 \right\}_{k \in \mathbb{Z}}
\right\|_{\ell^q} \right\|_{L^1}
\\&\lesssim     
\left\| \left\|
 \left\{
   M( |  \dot{\Delta}_{k} f_k |^{r_1   } )^{\frac{\gamma}{r_1} }
   M( | \dot{\Delta}_{k}   f_k |^{r_2  } )^{\frac{\delta}{r_2} }
 \right\}_{k \in \mathbb{Z}}
\right\|_{\ell^q} \right\|_{L^1}
\\&\lesssim    
\left\| \left\|
 \left\{
    M( |  \dot{\Delta}_{k} f_k |^{r_1   } )^{\frac{\gamma}{r_1} }
 \right\}_{k \in \mathbb{Z}}
\right\|_{\ell^{\frac{q}{\gamma} }} \right\|_{L^{ \frac{1}{\gamma}   }} \left\| \left\|
 \left\{
    M( | \dot{\Delta}_{k}   f_k |^{r_2  } )^{\frac{\delta}{r_2} }
 \right\}_{k \in \mathbb{Z}}
\right\|_{\ell^{\frac{q}{\delta} }} \right\|_{L^{ \frac{1}{\delta} }}
\\&=          
\left\| \left\|
 \left\{
    \dot{\Delta}_{k}   f_k
 \right\}_{k \in \mathbb{Z}}
\right\|_{\ell^{q }}^{\gamma}  \right\|_{L^{\frac{1}{\gamma} } }
\left\|
 \left\|
  \left\{
       \dot{\Delta}_{k} f_k
  \right\}_{k \in \mathbb{Z}}
 \right\|_{\ell^{q }}^{\delta}
\right\|_{L^{\frac{1}{\delta}  }}
\\&=           
\left\| \left\|
 \left\{
 \dot{\Delta}_{k}   f_k
 \right\}_{k \in \mathbb{Z}}
\right\|_{\ell^{q }} \right\|_{L^{1} },
\end{align*}
where
$\gamma + \delta =1 $,
 $\frac{\gamma}{r_1}  + \frac{\delta}{r_2} = \frac{1}{r}$.\footnote{In the following,
 $\gamma,\delta, r_1, r_2$ and $r$ are always positive real numbers.}
This implies the estimate.
\hfill$\Box$\par

\medskip

Equipped with Proposition \ref{prop:3} and its corollaries, we can
verify a Moser type inequality in $F_{1,\,\infty}^{s}(\mathbb{R}^d)$
and a commutator estimate in
$\dot{F}_{1,\infty}^{s}(\mathbb{R}^d)$. The ways of proofs have been
well-developed. The proofs are placed at Appendices.

\begin{proposition}[Moser type inequality in $F_{1,\infty}^{s}(\mathbb{R}^d)$]
\label{Moser}
Let $s > 0$. For scalar functions $f$ and $g$, we have
\begin{equation}
\label{eq:4}
\norm{fg}_{{F}_{1,\infty}^{s}}
\lesssim
    \norm{f}_{L^\infty} \norm{g}_{{F}_{1,\infty}^{s}}
    +
    \norm{g}_{L^\infty} \norm{f}_{{F}_{1,\infty}^{s}}.
\end{equation}
\end{proposition}

The bracket operator  $[\cdot, \cdot]$ below is defined by
$$
\left[
    f,\,\dot{\Delta}_{j}
\right] g := f
    \dot{\Delta}_{j}g
         -\dot{\Delta}_{j}( fg ).
$$
\begin{proposition}[Commutator estimate in $\dot{F}_{1,\infty}^{s}(\mathbb{R}^d)$]
\label{comm_est} Let $s > 0$. For a scalar function $f$ and a
divergence-free vector field $u$, we have
\begin{equation}
\label{eq:9}
\int_{\mathbb{R}^{d}}
\!
    \underset{j\in\mathbb{Z}}{\sup}
        |2^{js}
        ([u,\,\dot{\Delta}_{j}],  \nabla) f
        |(x)
dx \lesssim
    \norm{\nabla u}_{L^\infty}
    \norm{f}_{\dot{F}_{1,\infty}^{s}}
      \!\!  +\norm{u}_{\dot{F}_{1,\infty}^{s}}
        \norm{\nabla f}_{L^\infty}.
\end{equation}
\end{proposition}

\subsection{Boundedness of Leray projection operator in  $\dot{F}_{1,\infty}^{s}(\mathbb{R}^d)$}
\label{Leray}

We treat the continuity of  the Leray projection
$\mathbb{P}$ in $\dot{F}_{1,\infty}^{s}(\mathbb{R}^d)$,
which is defined by
$$
\mathbb{P}(u):= u-\nabla \Delta^{-1}(\nabla\cdot u),
$$
for an appropriate vector field $u$ in $\mathbb{R}^{d}$.
We present the following lemma:
\begin{lemma}\label{lem-Leray}
For any $k = 1, 2, \cdots, d$,  we have
\begin{align*}
\|
   \partial_{k} \Delta^{-1} u
\|_{\dot{F}_{1,\infty}^{s}} \lesssim \| u
\|_{\dot{F}_{1,\infty}^{s-1}}.
\end{align*}
\end{lemma}
\noindent {\bf Proof.}   
The corresponding symbol $m_k$ with respect to $\partial_{k}
\Delta^{-1}$ is $ m_{k}(\xi)=\frac{\xi_{k} }{|\xi|^2} $ with
 $\xi = (\xi_1, \xi_2, \cdots, \xi_d) \in\mathbb{R}^{d}$,
that is,
\begin{align}
\partial_{k} \Delta^{-1} u
=
\mathcal{F}^{-1}(m_{k}(\xi)\hat{u}). \label{lem3.2.5-1}
\end{align}
We choose a function
$\tilde{{\bf 1}} \in C_{0}^{\infty}(\mathbb{R}^{d})$ with
$\tilde{{\bf 1}}=1$ on $\mathrm{supp} \:
{h}_{0}=\{\xi:1/2\leq|\xi|\leq2\}$ and
$\tilde{{\bf 1}}(0)=0$.
From the fact that $m_k$ is a
symbol of homogeneous of order $-1$, we have
\begin{align}
\dot{\Delta}_{j} \mathcal{F}^{-1}
  (m_{k}(\xi)\,\hat{u})
&=
\mathcal{F}^{-1}(m_{k}(\xi))\ast\dot{\Delta}_{j}u  \label{exch}
\\&=
\mathcal{F}^{-1}
        (\tilde{{\bf 1}}(2^{-j}\xi)
            \times  2^{-j} m_{k}(2^{-j}  \xi))
    \ast\dot{\Delta}_{j}u
    \nonumber
\\&=
2^{-j} [\mathrm{K}]_{2^{j}} \ast  \dot{\Delta}_{j}u,  \nonumber
\end{align}
where $ \mathrm{K}:=\mathcal{F}^{-1}(\tilde{{\bf 1}}(\xi)\,m_{k}(
\xi)). $ For any positive integer $N_0$, we observe that
\begin{align*}
(1+|x|^{2})^{N_0} | \mathrm{K}(x)| &= \left|
\mathcal{F}^{-1}((I-\Delta)^{N_0}\tilde{{\bf 1}}(\xi)\,m_{k}(\xi))
\right|
\\&=
\left| \sum_{|\alpha|+|\beta|\leq 2N_{0}}c_{\alpha,\,\beta}
        \mathcal{F}^{-1}((\partial^{\alpha}\tilde{{\bf 1}})(\partial^{\beta}m_{k}))
\right| \lesssim 1,
\end{align*}
where $\alpha,\,\beta$ are multi-indices in $ \mathbb{R}^{d}$ and
$c_{\alpha,\,\beta}$ are integers whose exact values do not matter.
For $r > 0$, we may take an integer $N_0$ so large that $N_{0} >
\frac{d}{2} + \frac{d}{r}$ to have
\begin{align*}
 | \mathrm{K}(x) | ( 1+|x|)^{\frac{d}{r}}
\lesssim
\frac{(1+|x|)^{\frac{d}{r}}} {(1+|x|^{2})^{N_0}} \lesssim
\frac{1} {(1+|x|^{2})^{N_1}},
\end{align*}
where  $N_1$ is a suitable integer greater than $\frac{d}{2}$. Then the
function $\displaystyle\frac{1}{(1+|x|^{2})^{N_1}}$ is a nonnegative
radial decreasing integrable function, and hence
by Corollary \ref{coro-2} we derive that
\begin{align*}
\norm
    {
        \mathcal{F}^{-1}(m_{k}(\xi)\hat{u})
    }
_{\dot{F}_{1,\infty}^{s}}
\lesssim
\int
  _{\mathbb{R}^{d}}
\underset{j\in\mathbb{Z}}
  {\sup}
     \left|
         [\mathrm{K}]_{2^{j}} \ast (2^{j(s-1)} \dot{\Delta}_{j}u)
     \right|
(x)dx
\lesssim
\norm{u}_{\dot{F}_{1,\infty}^{s-1}}.
\end{align*}
From (\ref{lem3.2.5-1}), the arguments are completed.
\hfill$\Box$\par

\bigskip

Lemma \ref{lem-Leray} implies:
\begin{corollary} \label{Leray prop2}
The Leray projection $\mathbb{P}$ is continuous on
$\dot{F}_{1,\infty}^{s}(\mathbb{R}^d)$, that is,
\begin{align*}
\norm{ \mathbb{P} u }_{\dot{F}_{1,\infty}^{s}} \lesssim
\norm{u}_{\dot{F}_{1,\infty}^{s}}.
\end{align*}
\end{corollary}

\section{Persistence of the solution in $F_{1,\infty}^{s}(\mathbb{R}^d)$}
\label{persist_sol}

In this section we will prove the first and second properties
of the main theorem stated at page \pageref{Main Theorem}.
For $s \geq d +1$, a divergence free vector field
$$
u_0 \in  F_{1,\infty}^{s}(\mathbb{R}^d)
$$
is given.
According to the well-known result in \cite{P-P1} together with the fact that
$u_0 \in
F_{1,\infty}^{s}(\mathbb{R}^d)
 \subset
 B_{\infty, 1}^{1}(\mathbb{R}^d)$(\cite{Jaw} and Bernstein's lemma),
 there exists a time $T >0$ and a solution $u$
to the Euler equation (\ref{Euler}) in
$
C([0, T); B_{\infty, 1}^{1}(\mathbb{R}^d)  )
$.\footnote{   
The temporal continuity is due to the
denseness of smooth functions in $B_{\infty, 1}^{1}(\mathbb{R}^d)$.
See the arguments in \cite{P-P1} or \cite{P-P3}.
}    
Thanks to this,  it suffices to
prove that the solution $u(t)$ stays in
$F_{1,\infty}^{s}(\mathbb{R}^d)$ for a while.

First, we consider the trajectory map $X(x,\,t)$ defined by
\begin{align}
\frac{\partial}{\partial t} X(x,\,t)
 &= u(X(x,\,t),\,t),    \label{traj_flow}     \\
X(x,\,0)
&=x.
\nonumber
\end{align}
Then
from the Euler equation
(\ref{Euler}), we have
\begin{equation*}
u(X(x,\,t),\,t) = u(x)-\int_{0}^{t}\nabla p(X(x,\,\tau),\,\tau)
d\tau,
\end{equation*}
and so
\begin{align}
\norm{u}_{L^{1}}
\leq
\norm{u_{0}}_{L^{1}} +\int_{0}^{t}
    \norm{\nabla p}_{L^{1}} d\tau.  \label{L-1-est}
\end{align}
In order to estimate  the pressure term, we recall that
$$
\nabla p=\nabla(-\Delta)^{-1}\mathrm{div}((u,\,\nabla)u).
$$

\begin{lemma}\label{lem-pres}
For  $s > 0$, we have
\begin{align*}
\|
\nabla p
\|_{L^{1}}
\lesssim
\| u \|_{{ W}^{1, \infty}}
\|
u
 \|_{\dot{F}_{1,\infty}^{s}}.
\end{align*}
\end{lemma}
\noindent {\bf Proof.}   
Note that
\begin{align}
\|
\nabla p
\|_{L^{1}}
\lesssim
\|
\Delta_{-1} \nabla p
 \|_{L^1}
 +
\|
\nabla p
 \|_{\dot{F}_{1,\infty}^{s}}. \label{pres-est}
\end{align}
Successive applications of Lemma \ref{lem-Leray}, Proposition \ref{Moser} and Remark
\ref{rmk-2} imply
\begin{align}
\norm{\nabla p}_{\dot{F}_{1,\infty}^{s}}
\! \lesssim \!
\norm{\mathrm{div}(u,\,\nabla)u}_{\dot{F}_{1,\infty}^{s-1}}
\nonumber
&\! \lesssim
\sum_{i,\,j=1}^{d}
    \norm{\partial_{j} u^{i}
          \partial_{i} u^{j}}
    _{\dot{F}_{1,\infty}^{s-1}}
\nonumber
\\
&\!\lesssim
\norm{\nabla u}_{L^{\infty}} \!
\norm{\nabla
}_{\dot{F}_{1,\infty}^{s-1}}
\lesssim
\norm{\nabla u}_{L^{\infty}} \!
\| u \|_{\dot{F}_{1,\infty}^{s}}.\label{eest2}
\end{align}

To estimate the term $\norm{ \Delta_{-1} \nabla p  }_{{\bf
L}^{1}}$, we observe that
\begin{align*}
\norm{ \Delta_{-1}  \nabla p  }_{{ L}^{1} }
= \norm{  (\nabla \Phi) * p  \,
}_{{ L}^{1}}
= \norm{ (\Delta^{-1}\nabla\Phi) \ast
\mathrm{div}(u,\,\nabla)u }_{{ L}^{1}}.
\end{align*}
Applying  the fundamental solution of the Laplace's equation,
we get that for $d\geq 3$,
\begin{align}
\norm{ \Delta_{-1}  \nabla p  }_{{ L}^{1} }
&\leq
\frac{\Gamma(d/ 2 + 1)}{\sqrt{ \pi}^{d} d({d-2)}}
\sum_{i, j =1}^{d}
\bigg(
    \norm{
            \left(
                \tilde{{\bf 1}} \nabla\frac{1}{|x|^{d-2}}
            \right)
                \ast
                (\partial_{i} \partial_{j} \Phi)
           }_{{ L}^{1}}
\nonumber
\\
& \phantom{...........}
 + \norm{
        \partial_{i} \partial_{j}
            \left(
                (1- \tilde{{\bf 1}} )
                \nabla\frac{1}{|x|^{d-2}}
            \right)
        \ast\Phi
       }_{{ L}^{1}}
\bigg) \norm{u \otimes  u}_{{ L}^{1}}
\nonumber
\\
&\lesssim
\norm{u \otimes  u}_{{L}^{1}}
\lesssim
\| u \|_{{ L}^{\infty}}
\| u \|_{{ L}^{1}}
\lesssim
\| u \|_{{ L}^{\infty}}
\| u \|_{\dot{F}_{1,\infty}^{s}},  \label{eest2-2}
\end{align}
where $ \tilde{{\bf 1}} \in C_{0}^{\infty}(\mathbb{R}^{d})$ with
$\tilde{{\bf 1}}(x)= 1$  for  $|x|\leq1$.
We have
the same result for the two dimensional case $d= 2$.
Plugging the estimates (\ref{eest2}) and (\ref{eest2-2})
in the inequality (\ref{L-1-est}), we get the result.
\hfill$\Box$\par

\bigskip

Lemma \ref{lem-pres} leads to
\begin{align}\label{est5}
\norm{u}_{L^{1}} \lesssim \norm{u_{0}}_{L^{1}} +\int_{0}^{t}
  \| u \|_{{ W}^{1, \infty}}
   \|
    u
   \|_{\dot{F}_{1,\infty}^{s}}
   d\tau.
\end{align}
By applying $\dot{\Delta}_{j}$ and adding the term
$(u,\,\nabla)\dot{\Delta}_{j}u$ on both sides of the Euler equation
(\ref{Euler}) at $(X(x,\,t),\,t)$, we get
\begin{align}
\frac{d}{dt}\!\!\left[\dot{\Delta}_{j}u(X(x,\,t),\,t)\right]\!\!
&=-\nabla \dot{\Delta}_{j}p(X(x,\,t),t)
    +([u,\,\dot{\Delta}_{j}],\,\nabla)u(X(x,\,t), t).
\label{pse-1}
\end{align}
Integrate with respect to time $t$ over $[0, t]$,
 and then take $F^{s}_{1, \infty}$-norm
 on both sides of
(\ref{pse-1})
 to obtain
\begin{align}
\!\!\norm{u}_{\dot{F}_{1,\infty}^{s}}
 \!\!\leq \!
 \norm{u_{0}}_{\dot{F}_{1,\infty}^{s}}
\!\!\!\!+ \!\!\int_{0}^{t}\!\!
    \norm{\nabla p}_{\dot{F}_{1,\infty}^{s}} \!\!\!d\tau
\!+\!\!
\int_{0}^{t} \!\!\!
    \int_{\mathbb{R}^{d}}
        \underset{j\in\mathbb{Z}}
          {\sup}
          |2^{js}([u,\,\dot{\Delta}_{j}],\nabla)u|(x)dx
         d\tau. \label{est1}
\end{align}
Owing to the estimate (\ref{eest2}) and
 Proposition \ref{comm_est},
we can  derive that
\begin{align}
\label{est4}
\norm{u}_{\dot{F}_{1,\infty}^{s}} \lesssim
\norm{u_{0}}_{\dot{F}_{1,\infty}^{s}} +\int_{0}^{t} \norm{\nabla
u}_{L^{\infty}} \norm{ u}_{\dot{F}_{1,\infty}^{s}} d\tau.
\end{align}
Hence, combining (\ref{est4}) and (\ref{est5}), we find that
\begin{align}
\norm{u}_{F_{1,\infty}^{s}}
&\lesssim
\norm{u_{0}}_{F_{1,\infty}^{s}}
+
\int_{0}^{t}
\| u \|_{{ W}^{1, \infty}}
\norm{u}_{F_{1,\infty}^{s}} d\tau
\label{est-f}
\\
&\lesssim
\norm{u_{0}}_{F_{1,\infty}^{s}}
+
\int_{0}^{t}
\norm{u}_{F_{1,\infty}^{s}}^{2} d\tau.
\label{est-f-1}
\end{align}
 By virtue of  Gronwall's inequality, (\ref{est-f-1}) leads to
\begin{align}
\sup_{0 \leq \tau \leq t} \| u(\tau) \|_{F_{1,\infty}^{s}} \leq \:
 C_0 \| u_0 \|_{F_{1,\infty}^{s}}
 \exp\left\{C_0 \int_0^t \sup_{0 \leq \tau' \leq \tau}
       \|  u(\tau') \|_{ F_{1,\infty}^{s}} d \tau \right\}
       \label{blowup}
\end{align}
for some positive constant $C_0$.
Let
\[
y(t) :=
 C_0 \| u_0 \|_{F_{1,\infty}^{s}}
 \exp\left\{C_0 \int_0^t \sup_{0 \leq \tau' \leq \tau}
       \|  u (\tau') \|_{ F_{1,\infty}^{s} } d \tau \right\}.
\]
Then from (\ref{blowup}), we note that
\begin{align}
\frac{d}{d t} y(t)  \leq  C_0 y(t)^2, \qquad y(0) =  C_0  \| u_0
\|_{F_{1,\infty}^{s}}. \label{odi}
\end{align}
Solving the separable ordinary differential inequality (\ref{odi}),
we see  that
\[
y(t)
      \leq
\frac{C_0  \| u_0 \|_{F_{1,\infty}^{s}}} {1 - t C^2_0  \| u_0
\|_{F_{1,\infty}^{s}} },
\]
if
$0 \leq t
<
\frac{1}{C^2_0  \| u_0 \|_{F_{1,\infty}^{s}}  } := T_0$.
Hence it is shown that
\begin{align}
\sup_{0 \leq \tau \leq t} \| u(\tau) \|_{F_{1,\infty}^{s}}
 \leq \ y(t), \quad t \in [0, T_0 ).
              \label{est-bdd}
\end{align}
This illustrates that
 the solution $u(t)$ stays in $F_{1,\infty}^{s}$, at least, for $t \in [0, T_0)$.
 This completes the proof of the local existence and uniqueness of the solution.

\bigskip

For the two dimensional case $d = 2$,
 the estimate (\ref{est-f}) can be rewritten as
\begin{align}
 \| u(t) \|_{ {F}^{s}_{1, \infty}}
&\leq \:
 C \| u_0  \|_{ {F}^{s}_{1, \infty}}
 \exp
  \left\{C
    \int_0^t
       \| u (\tau) \|_{ {W}^{1, \infty}}
      d \tau
  \right\} \label{est-2D-1}
\end{align}
for some positive constant $C$.
Then from the well-known facts that
$\|  u  \|_{ {W}^{1, \infty}}
\lesssim
\|
 u
\|_{ {B}^{1}_{\infty, 1} }$
and
$\| u(t) \|_{ B^1_{\infty, 1 } } < \infty$
for all $t \in [0, \infty) $,
we establish that
 $u(t)$ stays inside the space $ {F}^{s}_{1, \infty}(\mathbb{R}^d)$
for all time $t$.
\hfill$\Box$\par

\section*{Appendices}
\label{Appendices}

We present the proofs of
Proposition \ref{Moser} and Proposition \ref{comm_est}. The key point is
to exhibit a detour way to avoid the discontinuity of the Hardy-Littlewood
 maximal operator
$M: L^1(\mathbb{R}^d) \to L^1(\mathbb{R}^d)$ with respect to $L^1$-norm.
The arguments are borrowed from \cite{Guo-Li}.

Bony's paraproduct formula
decomposes the product $f g$ into three parts as follows:
\begin{align}
fg = T_{f}g + T_{g}f
  + R\left(f,\,g\right),
\label{Bony}
\end{align}
where the para-product $T_{f}g$ and the remainder $R\left(f,\,g\right)$
 defined by
$$
T_{f}g
:=
 \sum_{j \in\mathbb{Z}}
 S_{j-4}f \Delta_{j}g
\quad \mbox{and} \quad
R \left(f,\,g\right) =
\sum_{\left|i-j\right|\leq3}\Delta_{i}f\Delta_{j}g,
$$
respectively \cite{B}.
Then we note  that
\begin{align}
\Delta_{k} T_{f}g &= \sum_{\left|j-k\right|\leq2}
     \Delta_{k}
       \left(S_{j-4}f\Delta_{j}g\right),
\label{para-prod}
\end{align}
because $\Delta_{k}\left(S_{j-4}f\Delta_{j}g\right)=0$
if $\left|j-k\right|\geq3$.
We also see
that
\begin{align}
\Delta_{k} R(f, g)
=
    \sum_{j= k-5}^{\infty}
 \sum_{\left|l\right|\leq3}
        \Delta_{k}
        \left(\Delta_{j}f\Delta_{j+l}g\right)
\label{remainder-1}
\end{align}
since
$
\mathrm{supp} \:\mathcal{F}
 \left(\Delta_{k}\left(\Delta_{j}f\Delta_{j+l}g\right)\right)
= \varnothing
$
if $ j \leq k-6$.

Similarly, the homogeneous version
of the Bony's paraproduct formula is of the
form:
\begin{align}
fg = \dot{T}_{f}g + \dot{T}_{g}f
  + \dot{R} \left(f,\,g\right),
\label{homoBony}
\end{align}
where the (homogeneous) Bony's para-product $T_{f}g$
and
(homogeneous) remainder
 defined by
$$
\dot{T}_{f}g :=
 \sum_{j \in\mathbb{Z}}
 \dot{S}_{j-4}f \dot{\Delta}_{j}g
\quad \mbox{and} \quad
\dot{R} \left(f,\,g\right) = \sum_{\left|i-j\right|\leq3}
  \dot{\Delta}_{i}f
  \dot{\Delta}_{j}g.
$$

\noindent {\bf Proof of Proposition \ref{Moser}.}
Bony's paraproduct decomposition can be read as
\[
fg = T_{f}g + T_{g}f
     + R\left(f,\,g\right).
\]
By the formula \eqref{para-prod} and
Corollary \ref{coro:1}, we have that for any $r\in(0,\,\infty)$,
\begin{align*}
    |2^{ks}(\Delta_{k}T_{f}g)|
&=\bigg|
    2^{ks}
    \sum_{\left|j-k\right|\leq2}
     \Delta_{k}
       \left(S_{j-4}f\Delta_{j}g\right)
\bigg|\\
&\lesssim \sum_{\left|l\right|\leq2}2^{ks} \left[
    2^{l\frac{d}{r}}M(|S_{k+l-4}f|)
    \{M(|\Delta_{k+l}g|^{r})\}^{\frac{1}{r}}
    (x)
\right].
\end{align*}
Then, for fixed $0<r<1$,
the continuity of vector-valued Hardy-Littlewood
maximal operator (the inequality \eqref{cont_max-fn}) can be employed to get
\begin{align}
\norm{T_{f}g}_{{F}_{1,\infty}^{s}} &=
 \int_{\mathbb{R}^{d}}\underset{k\in\mathbb{Z}}
 {\sup} \,
 2^{ks}
 \left|
   \Delta_{k}T_{f}g\right|(x)dx
\nonumber \\
&\lesssim
  \sum_{l=-2}^{2}
     \int_{\mathbb{R}^{d}}
      \underset{k\in\mathbb{Z}}
      {\sup}
\left[
        2^{2 \times \frac{d}{r}}
   M
   \left(
          \left|
             S_{k+l-4}f
          \right|
    \right)
\left\{
    M
    \left(
          \left|2^{ks}
             \Delta_{k+l}g
          \right|^r
    \right)
\right\}^{\frac{1}{r}}
    (x)
\right]
      dx
\nonumber \\
&\lesssim
 \norm{Mf}_{L^{\infty}}
\sum_{l=-2}^{2}
 \int_{\mathbb{R}^{d}}
    \underset{k \in \mathbb{Z}}{\sup} \;
     2^{ks}
     \left|\Delta_{k+l}g\right|
     \left(x\right)
 dx
\nonumber \\
&\lesssim \norm{f}_{L^{\infty}} \norm{g}_{{F}_{1,\infty}^{s}}.
\label{Moser:1}
\end{align}
Similarly, we can get
\begin{align}
\norm{T_{g}f}_{{F}_{1,\infty}^{s}} \lesssim \norm{g}_{L^{\infty}}
\norm{f}_{{F}_{1,\infty}^{s}}. \label{Moser:2}
\end{align}
For the remainder term,
the formula \eqref{remainder-1}
 and
 Corollary \ref{coro:1}
lead  to
\begin{align}
&\norm{R(f,\,g)}_{{F}_{1,\infty}^{s}}     \nonumber\\
& \!\!\leq \!\!\sum_{l=-3}^3 \!\!
    \int_{\mathbb{R}^d}
    \underset{k \in \mathbb{Z}}{\sup} \;
    \left|
        \sum_{j=k-5}^\infty
        2^{ks}
        \Delta_{k}(\Delta_{j}f \Delta_{j+l}g)
    \right|(x)dx
        \nonumber\\
& \!\!\lesssim \!\!\sum_{l=-3}^3 \!\!
    \int_{\mathbb{R}^d}
    \underset{k \in \mathbb{Z}}{\sup} \;
    \Bigg|
        \sum_{j=k-5}^\infty \!
        2^{(k-j-l)(s-\frac{d}{r})}
        M(\Delta_{j}f) \!
        \left[M(|2^{(j+l)s}\Delta_{j+l}g|^{r_1})\right]
            ^{\frac{\gamma}{r_1}}
        \nonumber\\
&\phantom{.................................................}
    \times
    \left[M(|2^{(j+l)s}\Delta_{j+l}g|^{r_2})\right]
            ^{\frac{\delta}{r_2}}    \Bigg|(x)dx,
            \label{Moser:3}
\end{align}
where $\gamma+\delta=1$ with
$\dfrac{\gamma}{r_1}+\dfrac{\delta}{r_2}=\dfrac{1}{r}$.
We choose
$r>0$
so large that we can have $s>\dfrac{d}{r}$.
The integrand of \eqref{Moser:3} is equal to
\[
\underset{k \in \mathbb{Z}}{\sup}
    \left| (a\ast b)(k) \right| (x),
\]
where the sequences $a:=\{a_{j}\}_{j\in\mathbb{Z}}$ and
$b:=\{b_{j}\}_{j\in\mathbb{Z}}$ are defined by
\[
    a_j:=\begin{cases}
        2^{(j-l)(s-\frac{d}{r})}, &\text{if } j\leq 5\\
        0, &\text{if } j>5
        \end{cases}
\]
and for $j\in\mathbb{Z}$,
\[
b_j:= M(\Delta_{j}f)
        \left[M(|2^{(j+l)s}\Delta_{j+l}g|^{r_1})\right]
            ^{\frac{\gamma}{r_1}}
        \times
        \left[M(|2^{(j+l)s}\Delta_{j+l}g|^{r_2})\right]
            ^{\frac{\delta}{r_2}}.
\]
Then, Young's inequality for $l^q$-sequences implies the
estimate
\begin{align*}
    \underset{k \in \mathbb{Z}}{\sup}
    \left| (a\ast b)(k) \right| (x)
    &\leq
        2^{-l(s-\frac{d}{r})}
        \left(\sum_{j=-\infty}^{5}
            2^{j(s-\frac{d}{r})}
        \right)
        \underset{j \in \mathbb{Z}}{\sup}
        |b_j|(x)
\lesssim
        \underset{j \in \mathbb{Z}}{\sup}
        |b_j|(x).
\end{align*}
Take $0<r_1,\,r_2<1$ and then H\"{o}ler's inequality and the
inequality \eqref{cont_max-fn} yield
\begin{align}
&\norm{R(f,\,g)}_{{F}_{1,\infty}^{s}}
\nonumber\\
&\lesssim \sum_{l=-3}^{3}
    \int_{\mathbb{R}^d}
    \underset{j \in \mathbb{Z}}{\sup} \;
    \Bigg|
        M(\Delta_{j}f)
        \left[M(|2^{(j+l)s}\Delta_{j+l}g|^{r_1})\right]
            ^{\frac{\gamma}{r_1}}
\nonumber\\
&\phantom{\lesssim \sum_{l=-3}^{3}
    \int_{\mathbb{R}^d}
    \underset{j \in \mathbb{Z}}{\sup} \;
    |M(\Delta_{j}f)|}
    \times
        \left[M(|2^{(j+l)s}\Delta_{j+l}g|^{r_2})\right]
            ^{\frac{\delta}{r_2}}
    \Bigg|(x)dx
\nonumber\\
&\lesssim \norm{Mf}_{L^{\infty}}\sum_{l=-3}^{3} \left(
    \int_{\mathbb{R}^d}
    \left|
        \underset{j \in \mathbb{Z}}{\sup}
        \left[M(|2^{(j+l)s}\Delta_{j+l}g|^{r_1})\right]
            ^{\frac{\gamma}{r_1}}
    \right|^{\frac{1}{\gamma}}
    (x)dx
\right)^{\gamma}
\nonumber\\
&\phantom{\lesssim\norm{Mf}_{L^{\infty}}\sum_{l=-3}^{3}}
    \times
    \left(
    \int_{\mathbb{R}^d}
    \left|
        \underset{j \in \mathbb{Z}}{\sup}
        \left[M(|2^{(j+l)s}\Delta_{j+l}g|^{r_2})\right]
            ^{\frac{\delta}{r_2}}
    \right|^{\frac{1}{\delta}}
    (x)dx
\right)^{\delta}
\nonumber\\
&\lesssim \norm{Mf}_{L^{\infty}}\sum_{l=-3}^{3}
    \int_{\mathbb{R}^d}
    \underset{j \in \mathbb{Z}}{\sup} \;
    |2^{(j+l)s}\Delta_{j+l}g|
    (x)dx
\lesssim
    \norm{f}_{L^{\infty}}
    \norm{g}_{{F}_{1,\,\infty}^{s}}.
    \label{Moser:4}
\end{align}
Combining the estimates \eqref{Moser:1}, \eqref{Moser:2} and
\eqref{Moser:4}, we obtain
the result (\ref{eq:4}).
\hfill$\Box$\par   

\bigskip

\noindent{\bf Proof of Proposition \ref{comm_est}.}
Let $u^l$ denote the $l$-th component of $u$ for
$1\leq l\leq d$.
By the homogeneous Bony's paraproduct decomposition, we
have
\begin{align*}
[u,\,\Dot{\Delta}_{j}]\cdot\nabla f
&=
    \sum_{l=1}^{d}
    \dot{T}_{\Dot{\Delta}_{j}\partial_{l}f} u^{l}
    +
    \sum_{l=1}^{d}
    \dot{R}(u^l,\,\Dot{\Delta}_{j}\partial_{l} f)
    +
    \sum_{l=1}^{d}
    [\dot{T}_{u^l},\,\Dot{\Delta}_{j}]\partial_{l} f
\\&\qquad\qquad\qquad
    -\sum_{l=1}^{d}
    \Dot{\Delta}_{j}
      \dot{T}_{\partial_{l}f} u^l
    -\sum_{l=1}^{d}
    \Dot{\Delta}_{j}
      \dot{R}(u^l,\,\partial_{l} f)
\\&
    :=\rom{1}+\rom{2}+\rom{3}+\rom{4}+\rom{5}.
\end{align*}

\underline{Estimate of $\rom{1}$} :
By the fact that
$
\dot{S}_{k-4}
 \Dot{\Delta}_{j} f = 0
$
if
$  k \leq j+2 $
and Young's inequality for $\ell^q$-series with
arguments used at page \pageref{Moser:3},
the first term  ($\rom{1}$)
can be estimated as
\begin{align*}
\int_{\mathbb{R}^{d}}
    \sup_{j\in\mathbb{Z}}
        |2^{js} \;  (\rom{1}) |(x)
dx
&\leq \sum_{l=1}^d \int_{\mathbb{R}^{d}}
    \sup_{j\in\mathbb{Z}}
        \left|
            \sum_{k=j+3}^{\infty}
                2^{js}
                (\dot{S}_{k-4}
                \Dot{\Delta}_{j}\partial_{l} f)
                (\Dot{\Delta}_{k}u^l)
        \right|(x)
  dx
\\
&\lesssim \sum_{l=1}^d
         \norm{\partial_{l} f }_{L^{\infty}}  \!\!
        \int_{\mathbb{R}^{d}}
            \sup_{j\in\mathbb{Z}}
                \left|
                    \sum_{k =j+3}^{\infty}
                        2^{(j-k)s}
                        (2^{ks}\Dot{\Delta}_{k}u^l)
                \right|(x)
        dx
\\&\lesssim
\sum_{l=1}^d
 \norm{\partial_{l} f }_{L^{\infty}}
\int_{\mathbb{R}^{d}}
    \left|
        \sum_{j= 3}^{\infty}
            2^{-js}
    \right|
    \sup_{j\in\mathbb{Z}}
        |2^{js}
        \Dot{\Delta}_{j}u^l
    |(x)
dx
\\&\lesssim
\norm{\nabla f}_{{\bf L}^{\infty}} \norm{u}_{\dot{\bf
F}_{1,\infty}^{s}}.
\end{align*}

\underline{Estimate of $\rom{2}$} :
We note that
 $ \mathrm{supp} \: \mathcal{F} (
    \Dot{\Delta}_{k}
        (
            \Dot{\Delta}_{j} f
        )
)
 \neq
  \varnothing $
 if
$k\leq j-2$ or $k\geq j+2$ to have
\begin{align*}
&\int_{\mathbb{R}^{d}}
    \underset{j\in\mathbb{Z}}{\sup}
        |2^{js}  (\rom{2}) |(x)
dx   \\
&\leq
 \sum_{l=1}^{d} \sum_{m=-3}^{3} \int_{\mathbb{R}^{d}}
    \sup_{j\in\mathbb{Z}}
        \left|
            2^{js}
            \sum_{k=j-1}^{j+1}
            \Dot{\Delta}_{k+m}u^l
            \Dot{\Delta}_{k}(\Dot{\Delta}_{j}\partial_{l}  f)
        \right|
(x)dx.
\end{align*}
Noting that
$
\mathrm{supp}\:
\mathcal{F} (
    (
        \Dot{\Delta}_{k+m}u^l
    )
    (
        \Dot{\Delta}_{k}
        (
            \Dot{\Delta}_{j}\partial_{l}f
        )
    )
)
\subseteq
\{
    \xi:|\xi|\leq2^{k+5}
\},
$
we can choose
$
\hat{\tilde{\mathbf{1}}} \in C_{0}^\infty
\left(
    \mathbb{R}^{d}
\right)
$
 with  $\hat{\tilde{\mathbf{1}}}= 1$ on $\{ \xi:|\xi|\leq2^{5}  \}$
 to get
\begin{align*}
    \Dot{\Delta}_{k+m}u^l
    \Dot{\Delta}_{k}
    (
        \Dot{\Delta}_{j}\partial_{l}f
    )
= [\tilde{\mathbf{1}}]_{2^{k}} \ast\partial_{l} (
    \Dot{\Delta}_{k+m}u^l
    \Dot{\Delta}_{k}
    (
        \Dot{\Delta}_{j}f
    )
)
\end{align*}
because $u$ is divergence free. Therefore, by
Corollary \ref{coro:1},
Bernstein lemma, Young's inequality and
Remark \ref{rmk-2}, we have that for $0 < r < 1 $,
\begin{align*}
&\int_{\mathbb{R}^{d}}
    \underset{j\in\mathbb{Z}}{\sup}
        |2^{js} \; (\rom{2}) |
(x)dx
\\
&\lesssim \sum_{l=1}^{d} \sum_{m=-3}^{3} \sum_{i=-1}^{1}
\int_{\mathbb{R}^{d}}
    \underset{j\in\mathbb{Z}}{\sup}
        \left|
            M
                (
                    |
                        2^{j}\Dot{\Delta}_{j+m+i}u^l
                    |)
             [
             M
                    (|
                        \Dot{\Delta}_{j+i}(2^{js}\Dot{\Delta}_{j}f)
                    |^{r}
                )
            ]^{\frac{1}{r}}
        \right|
(x)dx
\\
&\lesssim \norm{\nabla u}_{{\bf L}^\infty} \sum_{i=-1}^{1}
\int_{\mathbb{R}^{d}}
    \underset{j\in\mathbb{Z}}{\sup}
        \left|
                    |
                        \Dot{\Delta}_{j+i}(2^{js}\Dot{\Delta}_{j}f)
                    |^{r}
        \right|^{\frac{1}{r}}
(x)dx,   \\
&\lesssim
\norm{\nabla u}_{{\bf L}^\infty}
 \sum_{i=-1}^{1} \int_{\mathbb{R}^{d}}
    \underset{j\in\mathbb{Z}}{\sup}
        \left|
               M(     |
                        2^{js}\Dot{\Delta}_{j}f
                    |^{r}
                )
        \right|^{\frac{1}{r}}
(x)dx
\lesssim
\norm{\nabla u}_{{\bf L}^\infty}
\norm{f}_{\dot{F}_{1,\infty}^{s}}.
\end{align*}

\underline{Estimate of $\rom{3}$} :  
From the identity  (\ref{para-prod}) and the fact that
$$
 \mathrm{supp} \: \mathcal{F} (
    (\dot{S}_{k-4} u^l)
    (\Dot{\Delta}_{j}\Dot{\Delta}_{k}\partial_{l} f)
   )=\varnothing
$$
 if $k\geq j+2$ or $k\leq j-2$, we have
\begin{align*}
[
    \dot{T}_{u^{l}},\,\Dot{\Delta}_{j}
]
\partial_{l}f
    =\sum_{k=j-1}^{j+1}
        [
            \dot{S}_{k-4}u^{l},\,\Dot{\Delta}_{j}
        ]
            \Dot{\Delta}_{k}\partial_{l}f.
\end{align*}
Since the vector field $u$ is divergence free, we can derive
that
\begin{align*}
&\sum_{l=1}^{d} \sum_{k=j-1}^{j+1}
        [
            \dot{S}_{k-4}u^{l},\,\Dot{\Delta}_{j}
        ]
            \Dot{\Delta}_{k}\partial_{l} f
\\
&= - \!
\sum_{l=1}^{d} \!
\sum_{k=j-1}^{j+1} \!
\int_{\mathbb{R}^{d}}  \!\!\!
    2^{j(d+1)}
        (\partial_{l} \varphi_{0})
        (
            2^{j}(x-y)
        )
    (
        \dot{S}_{k-4}u^{l}(x)
        \! -  \!
        \dot{S}_{k-4}u^{l}(y)
    )
  \Dot{\Delta}_{k}
  f(y) dy.
\end{align*}
The mean value theorem provides
\begin{align*}
\dot{S}_{k-4}u^{l}(x)
        -\dot{S}_{k-4}u^{l}(y)
= (x- y) \cdot  (\nabla \dot{S}_{k-4}u^{l}) (\xi_{x}(y))
\end{align*}
for some $\xi_{x}(y)$ located on the line segment between $x$ and
$y$.
Applying Proposition \ref{prop:3}, we have
\begin{align*}
&\sum_{l=1}^{d}\sum_{k=j-1}^{j+1}
    \left|
        [
            \dot{S}_{k-4}u^{l},\,\Dot{\Delta}_{j}
        ]
            \Dot{\Delta}_{k}\partial_{l} f
    \right| (x)
\\
&=\sum_{l=1}^{d}
  \sum_{k=j-1}^{j+1}
        \left|
    \int_{\mathbb{R}^{d}}
        [\theta]_{2^j}(x-y)
        \cdot
        \left(
        \nabla \dot{S}_{k-4}u^{l} (\xi_{x} (y))
        \Dot{\Delta}_{k}f(y)
        \right)
     dy
        \right|
\\
&\lesssim
\sum_{m=-1}^{1}
    \sum_{l, i= 1}^d
     M \left(\partial_{i} \dot{S}_{k-4}u^{l} (\xi_{x} (\cdot) \right)
        M
            \left(
                    | \Dot{\Delta}_{j+m} f |^r
            \right)^{\frac{1}{r}}
            (x),
\end{align*}
where $0 < r < 1 $ and ${\bf \theta}(z) := z \partial_{l} \varphi(z)
\, \in {\mathcal S}(\mathbb{R}^d)$. Hence we conclude that
\begin{align*}
\int_{\mathbb{R}^{d}} \sup_{j\in\mathbb{Z}} &| 2^{js} \; (\rom{3})
|(x)
        dx
\\
&\lesssim \sum_{m=-1}^{1} \norm{\nabla u}_{{\bf
L}^{\infty}}
 \int_{\mathbb{R}^{d}}
    \underset{j\in\mathbb{Z}}{\sup}
        \left|
                M
                    \left(
                            |2^{js}\Dot{\Delta}_{j+m}f |^r
                    \right)^{\frac{1}{r}} (x)
        \right|
        dx
\\&\lesssim
\norm{\nabla u}_{{\bf L}^{\infty}}
\norm{f}_{\dot{F}_{1,\infty}^{s}}.
\end{align*}

\underline{Estimate of $\rom{4}$} :
The argument used at (\ref{Moser:1}) says that
 for $0 < r < 1 $,
\begin{align*}
&\int_{\mathbb{R}^{d}}
    \underset{j\in\mathbb{Z}}{\sup}
        |
            2^{js} \; (\rom{4})
        |
(x)dx
\lesssim
\norm{\nabla f}_{{\bf L}^{\infty}}
\norm{u}_{\dot{\bf F}_{1,\infty}^{s}}.
\end{align*}

\underline{Estimate of $\rom{5}$} :
The support of
$\mathcal{F} (
    \Dot{\Delta}_{j}
        (
            \Dot{\Delta}_{m}f
            \Dot{\Delta}_{m+k}g
        )
) $
is null set if $m<j-5$
as we pointed out at page \pageref{remainder-1}. So, by the
divergence-free condition of $\bf u$, we have that
\begin{align*}
(\rom{5})
 =
    -\sum_{l=1}^{d}
    \sum_{k=-3}^{3}
    \sum_{m=j-5}^{\infty}
    \Dot{\Delta}_{j}
    \partial_{l}
    (
        \Dot{\Delta}_{m}u^{l}
        \Dot{\Delta}_{m+k}f
    ).
\end{align*}
Applying Remark \ref{rmk-2} and Corollary \ref{coro:1}, we obtain that
\begin{align*}
&\int_{\mathbb{R}^{d}}
    \underset{j\in\mathbb{Z}}{\sup}
        |2^{js}  \; (\rom{5}) |
(x)dx
\\&\leq
 \sum_{l=1}^{d}  \sum_{k=-3}^{3}\int_{\mathbb{R}^{d}}
    \underset{j\in\mathbb{Z}}{\sup}
    \left|
         \sum_{m=j-5}^{\infty} 2^{(j-m-k)s}\Dot{\Delta}_{j}
       (
        2^{j}\Dot{\Delta}_{m}u^{l}
         \times 2^{(m+k)s}\Dot{\Delta}_{m+k}f
    )
    \right|(x)dx
\\&\lesssim
  \sum_{l=1}^{d}
  \sum_{k=-3}^{3}
  \int_{\mathbb{R}^{d}}
    \underset{j\in\mathbb{Z}}{\sup}
    \bigg|
    \sum_{m=j-5}^{\infty}
     2^{(j-m-k)(s-\frac{d}{r})}2^{j-m}
     \times
     M
    \left(
            2^{m}\Dot{\Delta}_{m}u^{l}
    \right)
\\&
  \hspace{0.4in}
  \times
\left[
     M
    \left(
            |
                2^{(m+k)s}\Dot{\Delta}_{m+k}f
            |
            ^{r_1}
    \right)
\right] ^{\frac{\gamma}{r_1}}
\left[
    M
        \left(
            |
                2^{(m+k)s}
                \Dot{\Delta}_{m+k}f
            |
            ^{r_2}
        \right)
\right] ^{\frac{\delta}{r_2}}
    \bigg|(x)dx,
\end{align*}
where
$\gamma + \delta =1 $  and
 $\frac{\gamma}{r_1}  + \frac{\delta}{r_2} = \frac{1}{r}$.
  We choose $r > 0$ so
large that we can have $s > \frac{d}{r}$,
and we also choose $0 < r_1, r_2 < 1$.
Then
 by Young's inequality
for $\ell^q$-series as the same argument at page \pageref{Moser:3}
together with  the Bernstein lemma,
 we get
\begin{align*}
\int_{\mathbb{R}^{d}}
    \underset{j\in\mathbb{Z}}{\sup}
        |2^{js} (\rom{5})|
(x)dx
&\lesssim
  \norm{\nabla u}_{{ L}^{\infty}}
   \sum_{k=-3}^{3}
    \left(
     \int_{\mathbb{R}^{d}}
      \underset{j\in\mathbb{Z}}{\sup}
        \left|
                         2^{js}\Dot{\Delta}_{j+k}f
        \right|
       (x)dx
     \right)
\\&\lesssim \norm{\nabla u}_{L^{\infty}}
\norm{f}_{\dot{F}_{1,\infty}^{s}}.
\end{align*}
Collecting the estimates of the terms (\rom{1}) $\sim$ (\rom{5})
altogether, we obtain the inequality (\ref{eq:9}).
\hfill$\Box$\par   

\end{document}